\documentclass[11pt,fleqn,twoside]{article}

\usepackage{amsmath}
\usepackage{amsthm}
\usepackage{amssymb}
\makeatletter
\newcommand{\prava}[1]{\small\it
\begin{flushleft}
Copyright \copyright \ 2000 by  #1
\end{flushleft}}

\newcommand{\name}[1]{\begin{flushleft}
                       \LARGE \bf #1
                       \end{flushleft}\vspace{-3mm}}

\newcommand{\Author}[1]{\begin{flushleft}
                       \it #1 \end{flushleft}}

\newcommand{\Adress}[1]{\begin{flushleft}
                       \it #1 \end{flushleft}}

\newcommand{\Date}[1]{\begin{flushleft}
                      \small  \it #1 \end{flushleft}}

\newcommand{\ehkol}{Author \ name}
\newcommand{\ohkol}{Article \ name}

\renewcommand{\@evenhead}{
\hspace*{-3pt}\raisebox{-15pt}[\headheight][0pt]{\vbox{\hbox to \textwidth 
{\thepage \hfil \ehkol}\vskip4pt \hrule}}}
\renewcommand{\@oddhead}{
\hspace*{-3pt}\raisebox{-15pt}[\headheight][0pt]{\vbox{\hbox to \textwidth 
{\ohkol \hfil \thepage}\vskip4pt\hrule}}}
\renewcommand{\@evenfoot}{}
\renewcommand{\@oddfoot}{}

\long\def\@makecaption#1#2{%
  \vskip\abovecaptionskip
  \sbox\@tempboxa{\small \textbf{#1.}\ \ #2}%
  \ifdim \wd\@tempboxa >\hsize
    {\small \textbf{#1.}\ \ #2}\par
  \else
    \global \@minipagefalse
    \hb@xt@\hsize{\hfil\box\@tempboxa\hfil}%
  \fi
  \vskip\belowcaptionskip}

\def\numberwithin#1#2{\@ifundefined{c@#1}{\@nocounterr{#1}}{%
  \@ifundefined{c@#2}{\@nocnterr{#2}}{%
  \@addtoreset{#1}{#2}%
  \toks@\@xp\@xp\@xp{\csname the#1\endcsname}%
  \@xp\xdef\csname the#1\endcsname
    {\@xp\@nx\csname the#2\endcsname
     .\the\toks@}}}}

\newcommand{\resetfootnoterule} {
  \renewcommand\footnoterule{%
  \kern-3\p@
  \hrule\@width.4\columnwidth
  \kern2.6\p@}
}

\makeatother

\numberwithin{equation}{section}

\def\pmq{\boldsymbol{q}}
\def\dps{\displaystyle}
\def\bn{{\bf{N}}}
\def\Z{{\bf{Z}}}

\makeatletter
\DeclareRobustCommand{\primfrac}[1]{%
  \PackageWarning{amsmath}{%
Foreign command \@backslashchar#1; %
\protect\frac\space or \protect\genfrac\space should be used instead%
  }
  \global\@xp\let\csname#1\@xp\endcsname\csname @@#1\endcsname
  \csname#1\endcsname
}
\makeatother

\begin{document}

\thispagestyle{empty}
\renewcommand{\ehkol}{B.A.\ Kupershmidt}
\renewcommand{\ohkol}{$q$-Probability:
I.\ Basic Discrete Distributions}

\begin{flushleft}
\footnotesize \sf
Journal of Nonlinear Mathematical Physics \qquad 2000, V.7, N~1,
\pageref{kuper_fp}--\pageref{kuper_lp}.
\hfill {\sc Article}
\end{flushleft}

\vspace{-5mm}

\renewcommand{\footnoterule}{}
{\renewcommand{\thefootnote}{}
 \footnotetext{\prava{B.A.\ Kupershmidt}}}

\name{$\pmq$-Probability: 
I.\ Basic Discrete Distributions}\label{kuper_fp}

\Author{Boris A.\ KUPERSHMIDT}

\Adress{The University of Tennessee Space Institute, 
Tullahoma, TN  37388, USA \\
e-mail:  bkupersh@utsi.edu}

%\Date{Received July 16, 1999; Revised September 11, 1999; Accepted
%September 16, 1999}
\Date{Received November 5, 1999; Accepted
December 9, 1999}

\begin{abstract}
\noindent
For basic discrete probability distributions, $-$ Bernoulli, 
Pascal, Poisson, hypergeometric, contagious, and uniform, $-$ 
$q$-analogs are proposed.
\end{abstract}

\section {Introduction}

$q$-analogs of classical formulae go back to Euler, $q$-binomial coefficients
were defined  by Gauss, and $q$-hypergeometric series were found by E. Heine in
1846.  The $q$-analysis was  developed by F. Jackson at the beginning of the
20$^{th}$ century, and the modern point of  view subsumes most of the old
developments into the subjects of Quantum Groups and Combinational  Enumeration.

The general philosophy of $q$-analogs is that of a deformation, with the
deformation  parameter $q$ being thought of as close to 1.  This point of view
is certainly not all-encompassing;  for example, representations of Quantum
groups when $q$ is a root of unity are of independent  interests; more
importantly, the $q$-pictures sometimes possess properties {\it{singular}} in 
($q-1$) or otherwise not regularly dependent on $(q-1)$; regularization of
divergent/infinite 
$(q=1)-$quantities is another useful feature of $q$-analogs... the list goes on.

The typical example is
\begin{equation}
  \lim_{x \to \infty} [x] = {1 \over 1 - q}, \ \ \ |q| < 1, \tag{1.1} 
\end{equation}
where
\begin{equation}
  [x] = [x]_q = {q^x - 1 \over q -1} \tag{1.2} 
\end{equation}
is the $q$-analog of a number (or object) $x$.  (A quick introduction to the
$q$-calculus is  available from many sources, e.g.\ Chapter 2 in [6].)

More examples of such sort will be found below in this paper, the 1$^{st}$ one
in a series  devoted to a $q$-probabilities.  In the next 6 sections we look at
$q$-analogs of Bernoulli,  Pascal, Poisson, hypergeometric, contagious, and
uniform distributions,  respectively.  It is surprising how many new 
effects appear compared to the classical theory at $q=1$.  For example, even for the Bernoulli 
distribution, the profound classical differences between finite and infinite
number of  trials are mitigated when $q$ enters the picture, so that one can
write down the probability  (formerly zero) of {\it{many}} individual events of
infinite type, such as
\begin{equation}
  (1 -^{\hskip-.076truein \cdot} p)^
  \infty = \ {{\dps{\sqcap^\infty_{i=0}}}} \  (1 - pq^i), \ \ \ |q|<1, 
  \tag{1.3} 
\end{equation}
the probability of coming up with all ``tails'' during an infinite number of
coin flips;  the probability of coming all tails during $n$ coin flips is 
\begin{equation}
  (1 -^{\hskip-.076truein \cdot} p)^n =  \ {{\dps{\sqcap^{n-1}_{i=0}}}}  \ (1 - pq^i). \tag{1.4} 
\end{equation}
More generally, the probability of coming up with precisely $\kappa$ heads
during an  infinite number of coin flips is
\begin{equation}
  {p^\kappa \over (1-q)...(1-q^\kappa)} (1-^{\hskip-.076truein \cdot} 
  p)^\infty, \ \ \ \kappa \in \bn. \tag{1.5} 
\end{equation}

The six probability distributions discussed in this paper are all discrete, 
univariate, basic, and relatively simple.  More classical distributions can be found in
[1, 3-5].

\section{$\pmq$-binomial distributions}

Suppose we have a random variable $\zeta \ -$ ``2-sided coin'' $-$ which takes two values:  1 \ with probability 
$p$, and $0$ with probability 
\begin{equation}
  p^\prime = 1 - p. \tag{2.1} 
\end{equation}
After $n$ throws, the total sum accumulated,
\begin{equation}
  \xi_n = \zeta_1 + ... + \zeta_n, \tag{2.2} 
\end{equation}
obeys the Bernoulli distribution
\begin{equation}
  Pr (\xi_n = \kappa) = {n \choose \kappa} p^\kappa p^{\prime n - \kappa}, \ \ \ 0 \leq \kappa 
  \leq n. \tag{2.3} 
\end{equation}

As a $q$-analog of this distribution we take (with $0 < q < 1)$
\begin{equation}
  Pr (\bar \xi_n = [\kappa]) = \bigg[{n \atop \kappa}\bigg] p^\kappa  
  (1 -^{\hskip-.09truein \cdot} p)^{n-\kappa} , \ \ \ 0 \leq \kappa \leq n,
  \tag{2.4}
\end{equation}
where
\begin{gather}
  \bigg[{x \atop \kappa}\bigg] = \bigg[{x \atop \kappa}\bigg]_q = {[x] ... 
  [x - \kappa + 1] \over [\kappa]!} , 
  \ \ \ \kappa \in \bn, \tag{2.5a} \\
  \bigg[{x \atop 0}\bigg] = 1, \ \ \ \bigg[{x \atop s}\bigg] = 0, 
  \ \ \ s \ \bar \in \ \Z_+, \tag{2.5b} 
\end{gather}
are the $q$-binomial coefficients, and
\begin{equation}
  [0]! = 1; \ \ [\kappa]! = [1] ... [\kappa], \ \ \ \kappa \in \bn, \tag{2.6} 
\end{equation}
are the $q$-factorials. 

To justify formula (2.4), we need to prove that
\begin{equation}
  \sum^n_{k =0} \bigg[{n \atop \kappa}\bigg] p^\kappa
(1 -^{\hskip-.076truein \cdot} 
  p)^{n-\kappa} = 1, \ \ \ 
  \forall p. \tag{2.7} 
\end{equation}
This formula follows from the following identity:
\begin{equation}
  \sum^n_{\kappa =0} \bigg[{n \atop \kappa}\bigg] p^\kappa
(1 -^{\hskip-.076truein \cdot} 
  v)^{n - \kappa} 
  = \sum^n_{\kappa =0} \bigg[{n \atop \kappa}\bigg] (p 
  -^{\hskip-.076truein \cdot} v)^{n-\kappa}, \tag{2.8} 
\end{equation}
when $v =p$; here
\begin{equation}
  (a \dot + b)^n : = \ \sqcap^{n-1}_{i=0} \ (a + q^i b), \ \ \ n \in \bn, \ \ (a \dot + b)^0 
  : = 1. \tag{2.9} 
\end{equation}
Formula (2.8) is, in turn, the $b = 1$-case of the general formula
\begin{equation}
  \sum^n_{\kappa = 0} \bigg[{n \atop \kappa}\bigg] a^\kappa (b \dot + v)^{n - \kappa} = \sum
  ^n_{\kappa = 0} \bigg[{n \atop \kappa}\bigg] b^\kappa (a \dot + v)^{n -
\kappa}. \tag{2.10} 
\end{equation}
To prove formula (2.10), let us use the easily checked by induction on $m$
Euler's formula
\begin{equation}
  (x \dot + y)^m = \sum^m_{j = 0} \bigg[{m \atop j}\bigg] x^{m-j} y^j q^{j \choose 2}. 
  \tag{2.11} 
\end{equation}
Then the LHS of (2.10) can be rewritten as
\begin{equation}
  \sum_{\kappa, \ell} \bigg[{n \atop \kappa}\bigg] a^\kappa \bigg[{n - \kappa \atop \ell}\bigg] 
  b^{n - \kappa - \ell} v^\ell q^{{\ell \choose 2}}, \tag{2.12L} 
\end{equation}
while the RHS of (2.10) can be similarly rewritten as
\begin{equation}
  \sum_{s, \ell} \bigg[{n \atop s}\bigg] b^s \bigg[{n - s \atop \ell}\bigg] a^{n-s-\ell} v^\ell 
  q^{\ell \choose 2}, \tag{2.12R}
\end{equation}
and these two double sums bijectively coincide, for each fixed $\ell$, when
$s$ is identified with 
$n - \ell - \kappa$.

Now, to calculate the expectation values of powers of $\bar \xi_n$, we use the
easily proved  by induction on $m \in \bn$ formula
\begin{equation}
  \bigg(x {d \over d_q x} \bigg)^m = \sum^m_{\kappa =1} {1 \over [ \kappa - 1]!} 
  \bigg(\sum^{\kappa -1}_{s=0} \bigg[{\kappa - 1 \atop s} \bigg] (-1)^s q^{s
  \choose 2} [\kappa  - s]^{m-1} \bigg) x^\kappa \bigg({d \over d_q x}
  \bigg)^\kappa, \tag{2.13} 
\end{equation}
where
\begin{equation}
  {df \over d_q x} : = {f(qx) - f(x) \over (q-1)x} \tag{2.14} 
\end{equation}
is the $q$-derivative:
\begin{equation}
  {d \over d_q x} (x^s) = [s] x^{s-1}. \tag{2.15} 
\end{equation}
In particular, for $m=2$, we get
\begin{equation}
  \bigg(x {d \over d_q x} \bigg)^2 = x {d \over d_q x} + q x^2 \bigg({d \over d_q x}\bigg)^2. 
  \tag{2.16} 
\end{equation}
Applying the operator ${\dps{\bigg(p {d \over d_q p} \bigg)^s\bigg|_{v = p}}},
\ s = 1,2,$ to formula (2.8), we obtain
\begin{gather}
  < \bar \xi_n > = E (\bar \xi_n) = \sum^n_{k=0} [\kappa] \bigg[{n \atop \kappa}\bigg] p^\kappa 
  (1 -^{\hskip-.076truein \cdot} 
  p)^{n-\kappa} = [n] p, \tag{2.17} \\
  < \bar \xi^2_n > = E(\bar \xi^2_n) = \sum^n_{\kappa =0} [\kappa]^2 \bigg[{n
  \atop \kappa}\bigg] 
  p^\kappa (1 -^{\hskip-.076truein \cdot} p)^{n-\kappa} = 
  [n] p+ q p^2 [n] [n-1], \tag{2.18} 
\end{gather}
where we used the obvious relation
\begin{equation}
  {d \over d_q p} (p \dot + v)^\ell = [ \ell] (p \dot + v)^{\ell-1}. \tag{2.19} 
\end{equation}
From formulae (2.17) and (2.18) we find that
\begin{equation}
  Var (\bar \xi_n) = < \bar \xi^2_n> - <\bar \xi_n>^2 \ = \ [n] p (1-p). \tag{2.20} 
\end{equation}

Notice that formulae (2.4), (2.17), (2.18), (2.20) have a well-defined limit
as $n \rightarrow \infty$:
\begin{gather}
  Pr (\bar \xi_\infty = [\kappa]) = {1 \over [\kappa]!} \bigg({p \over 1 - q} \bigg)^\kappa 
  (1 -^{\hskip-.076truein \cdot} p)^\infty, \tag{2.21} \\
  < \bar \xi_\infty > = {p \over 1 - q}, \tag{2.22} \\
  < \bar \xi^2_\infty> = {p \over 1 - q} + q \bigg({p \over 1 - q} \bigg)^2,
  \tag{2.23}  \\
  V ar (\bar \xi_\infty) = {p(1-p) \over 1-q}. \tag{2.24} 
\end{gather}

So far, we have treated the random variable $\bar \xi_n$ as an object in its
own right.  Let  us now turn to the representation of $\bar \xi_n$ as a sum of
$n$ ``coin'' throws, as reflected  in the classical formula (2.2).  This
sum-formula (2.2) remains intact under $q$-deformation.   However, for general
$q$, the random variables $\bar \zeta_i$'s are {\it{no \ longer \
independent}} or {\it{identically \ distributed}}.

More precisely, let
\begin{equation}
  \bar \zeta_1 = \begin{cases}1, & \text{with probability } p, \\
  0,& \text{with probability } 1 - p . \end{cases} \tag{2.25} 
\end{equation}
For $\kappa \in \bn$, suppose the random variables $\bar \zeta_1, ... , \bar
\zeta_\kappa$ have  already been defined.  Denote by $\alpha_\kappa =
\alpha_\kappa (\bar \zeta_1, ..., 
\bar \zeta_\kappa)$ the 
number of zeroes appearing among the values of the random variables $\bar \zeta_1, ..., 
\bar \zeta_\kappa$.  The random variable $\bar \zeta_{\kappa +1}$ takes the values $0, q^0, 
... , q^\kappa$, with the conditional probabilities 
\begin{gather}
  Pr (\bar \zeta_{\kappa +1}= 0 | \alpha_\kappa = r) = 1 - q^r p, \ \ \ 0 \leq r \leq \kappa, 
  \tag{2.26a} \\
  Pr (\bar \zeta_{\kappa+1} = q^\ell | \alpha_\kappa = r) =
  \delta^\ell_{\kappa-r} q^r p, \ \ \ 
  0 \leq r \leq \kappa. \tag{2.26b} 
\end{gather}
For example,
\begin{gather}
  Pr (\bar \zeta_2 = 0 | \bar \zeta_1 = 0) = 1- qp, \tag{2.27a} \\
  Pr (\bar \zeta_2 = 0 | \bar \zeta_1 = 1) = 1- p, \tag{2.27b} \\
  Pr (\bar \zeta_2 = 1 | \bar \zeta_1 = 0) = qp, \tag{2.27c} \\
  Pr (\bar \zeta_2 = q | \bar \zeta_1 = 1) = p. \tag{2.27d} 
\end{gather}
By induction on $\kappa$, it is easily seen that if the last before $\bar
\zeta_{\kappa+1}$  non-zero value appearing among $\zeta_1, ... \bar
\zeta_\kappa$  was $q^\ell$, then $\bar \zeta_{\kappa +1}$ can take only the
values $0$ and 
$q^{\ell+1}$ with non-zero probabilities; if all the $\bar \zeta_1, ... , \bar \zeta_\kappa$ took 
value $0$, then $\bar \zeta_{k+1}$ takes only the values $0$ and $1$ with non-zero probabilities.

A better description of the same distribution is possible, if instead of
conditional probabilities  we work with joint ones.  Denote by ${\bf{0}}^r$
the event of $r$ in a row appearances of zeroes, 
$r \in \Z_+$, and similarly by 
\begin{equation}
  {\bf{0}}^{a(0)} q^0 {\bf{0}}^{a(1)} q^1 ... {\bf{0}}^{a(k-1)} q^{\kappa-1} 
  {\bf{0}}^{a(\kappa)}, \ \ a (\cdot) \in 
  \Z_+, \tag{2.28} 
\end{equation}
the event of $a(0)$ of zeroes followed by $q^0=1$ followed by $a(1)$ zeroes
... Now set
\begin{equation}
  Pr ({\bf{0}}^{a(0)} q^0 ... q^{\kappa -1} {\bf{0}}^{a(\kappa)} ): = 
  (1 -^{\hskip-.076truein \cdot} p)^{\sum^{\kappa}_{0} a(s)} \
  \sqcap^{\kappa-1}_{i=1} (p q ^{\sum^{i}_{0} a(s)} ) ; \tag{2.29} 
\end{equation}
for $\kappa=0$, formula (2.29) is to be understood as
\begin{equation}
  Pr ({\bf{0}}^a) = (1 -^{\hskip-.076truein \cdot}  p)^a, \ \ \ a \in \bn.
  \tag{2.29$^\prime$} 
\end{equation}

Let us now verify that the ``microscopic'' formulae (2.29) imply the
``macroscopic'' formula  (2.4).  Denote
\begin{equation}
  |a (i) | = \sum^i_{s=0} a(s). \tag{2.30} 
\end{equation}
We have to verify that
\begin{equation}
  \sum_{|a(\kappa)| = n-\kappa} (1 -^{\hskip-.076truein \cdot} p)^{n-\kappa}
  p^\kappa \ 
  \sqcap^{\kappa-1}_{i=0} q^{|a(i)|} = \bigg[{n \atop \kappa}\bigg] p^\kappa 
  (1 -^{\hskip-.076truein \cdot}  p)^{n-\kappa}, \tag{2.31} 
\end{equation}
which is equivalent to the $q$-counting formula
\begin{equation}
  \sum_{|a(\kappa)| =n-\kappa} \ \sqcap^{\kappa-1}_{i=0} q^{|a(i)|} = \bigg[{n \atop \kappa}
  \bigg]. \tag{2.32} 
\end{equation}
(For $q=1$, we recover the classical result:  the number of solutions in
nonnegative integers  of the equation $a(0)+...+a(\kappa) = n - \kappa$ is
${\dps{{n \choose \kappa}}}$. )

We shall prove formula (2.32) by the double induction on $N:=n-\kappa$ and
$\kappa$,  in the form 
\begin{equation}
  \sum_{|a(\kappa)|=N} \ \sqcap^{\kappa-1}_{i=0} q^{|a(i)|} = \bigg[{N+\kappa \atop \kappa}\bigg].
  \tag{2.33} 
\end{equation}
Now, 
\begin{align}
  \sqcap^{\kappa-1}_{i=0} q^{|a(i)|} &= q^{a (0)} q^{a(0)+a(1) |} ... q^{a(0)
  + ... + 
  a(\kappa-1)} = q^{\sum^{k}_{0}(k-s)a(s)} \tag{2.34a} \\
  &= q^{\kappa|a(\kappa)|} q^{-\sum^{k}_{0} sa(s)}. \tag{2.34b} 
\end{align}
Therefore, the identity (2.33) becomes
\begin{equation}
  \sum_{|a(\kappa)|=N} q^{-\sum^{k}_{0}sa(s)} = q^{-Nk} \bigg[{N + \kappa \atop \kappa} \bigg]. 
  \tag{2.35} 
\end{equation}
For $N=0$, formula (2.35) becomes $1=1$ no matter what $\kappa $ is.  
For $\kappa = 0$, formula (2.35) is true by definition (2.29$^\prime)$; for $\kappa =1$, 
formula (2.35) becomes
\begin{equation*}
  \sum_{a(0) + a(1) = N} q^{-a(1)} = q^{-N} \bigg[{N + 1 \atop 1}\bigg], 
\end{equation*}
which is obviously true for all $N$.  Suppose formula (2.35) is true for the
pairs $(\kappa,  N=n)$ and $(\kappa - 1, N = n+1)$.  Consider the pair
$(\kappa, N = n+1)$.   Let's divide the set of the 
$a$'s with $|a (\kappa)|= n+1$ into two groups:  those with $a(\kappa)>0$ and those with 
$a (\kappa)=0$.  For the $1^{st}$ group,  the set $\bar a(0) = a(0), ... , \bar a(\kappa -1) 
= a(\kappa -1), 
\bar a(\kappa) = a (\kappa)-1$, satisfies \ $|\bar a (\kappa)|=n$, so that, by the induction 
assumption, 
\begin{equation}
  \sum q^{- \sum^{\kappa}_{0} sa(s)} = \sum_{|\bar a(\kappa)|=n} q ^{-\sum^{\kappa}_{0}
  s \bar a 
  (s)-\kappa} = q^{-\kappa} q^{- n \kappa} \bigg[{n + \kappa \atop
  \kappa}\bigg]. \tag{2.36} 
\end{equation}
The $2^{nd}$ group has effectively $\kappa-1 \ a$'s, so that again, by the
induction  assumption, 
\begin{equation}
  \sum q^{- \sum^{\kappa-1}_{0} sa(s)} = q^{-(n+1)(\kappa-1)} \bigg[{n+\kappa \atop \kappa -1}
  \bigg]. \tag{2.37} 
\end{equation}
We thus have to check that 
\begin{equation}
  q^{-\kappa (n+1)} \bigg[{n +\kappa \atop \kappa}\bigg] +q ^{-(n+1) (\kappa -1)}  
  \bigg[{n+\kappa \atop \kappa-1}\bigg] = q^{-(n+1)\kappa} \bigg[{n+1+\kappa
  \atop \kappa}\bigg], 
  \tag{2.38} 
\end{equation}
which is equivalent to
\begin{equation}
  \bigg[{n + \kappa \atop \kappa}\bigg] + \bigg[{ n + \kappa \atop \kappa - 1} \bigg] 
  q^{n+1} = \bigg[{n + 1 + \kappa \atop \kappa}\bigg], \tag{2.39} 
\end{equation}
which is obviously true.

Notice that for $a(\kappa)= \infty$, formulae (2.29) become
\begin{gather}
  Pr ({\bf{0}}^{a(0)} q^0 ... q^{\kappa-1} {\bf{0}}^\infty) = p^\kappa 
  (1 -^{\hskip-.076truein \cdot} p)^\infty q^{\sum^{\kappa-1}_{0}
  (\kappa-s)a(s)}, 
  \tag{2.40} \\
  Pr ({\bf{0}}^\infty) = (1 -^{\hskip-.076truein \cdot}  p)^\infty. 
  \tag{2.40$^\prime$}
\end{gather}

As in the classical theory (cf [10] p. 59), we can calculate the probability
of observing  
$\leq \kappa \ (<n)$ zeroes in $n$ trials:
\begin{align}
  Pr (\alpha_n \leq \kappa) &= \sum^\kappa_{i=0} \bigg[{n \atop i}\bigg]
  p^{n-i} (1 -^{\hskip-.076truein \cdot}  p)^i= [n] \bigg[{n-1 \atop
  \kappa}\bigg] \int^p_0 x^{n-1-\kappa} (1 -^{\hskip-.076truein \cdot} 
  qx)^\kappa d_q x \notag\\
  &= (\int^p_0 x^{n-1-\kappa} (1 -^{\hskip-.076truein
  \cdot} qx)^\kappa  d_q x) / (\int^1_0  x^{n-1-\kappa}
(1 -^{\hskip-.076truein
  \cdot} q x)^\kappa d_q x).
  \tag{2.41} 
\end{align}
Similarly, the probability of $\leq \ell \ (<n)$ non-zeroes in $n$ trials is
\begin{equation}
  Pr (\alpha_n > n-\ell) = \sum^\ell_{i = 0} \bigg[{n \atop i}\bigg] p^i 
  (1 -^{\hskip-.076truein \cdot} p)^{n-i} = 1 - [n] \bigg[{n-1 \atop \ell}
  \bigg] \int^p_0  x^\ell (1 -^{\hskip-.076truein \cdot} qx)^{n-1-\ell } d_q  x.
  \tag{2.42} 
\end{equation}
(Formulae (2.41) and (2.42) can be easily proven upon $q$-differentiation with
respect to 
$p$ and using the obvious relation
\begin{equation}
  {d \over d_q p} (1 -^{\hskip-.076truein \cdot} p)^n = - [n]
(1 -^{\hskip-.076truein \cdot}  qp
  )^{n-1} . \tag{2.43} 
\end{equation}
In the limit $n \rightarrow \infty$, this probability becomes
\begin{equation}
  \sum^\ell_{i=0} {1 \over [i]!} \bigg({p \over 1-q} \bigg)^i
(1 -^{\hskip-.076truein \cdot}  p)
  ^\infty =  1 - {1 \over 1 - q} \int^p_0 {1 \over [\ell]!} \bigg({x \over 1 -
  q}\bigg) ^\ell  (1 -^{\hskip-.076truein \cdot}  qx)^\infty d_q x. 
  \tag{2.44} 
\end{equation}

Many, if not all, classical formulae in probability have $q$-analogs.  Let's
take a look at a few  of such formulae involving higher moments for the
Bernoulli distribution.

First, applying the operator ${\dps{p^r \bigg({d \over d_qp}\bigg)^r \bigg|_{v
= p}}}$ to formula  (2.8) and using the relation
\begin{equation}
  x = [\kappa] \Rightarrow [\kappa -i] = q^{-i} (x - [i]), \tag{2.45} 
\end{equation}
we get
\begin{equation}
  E ( \sqcap^{\kappa-1}_{i=0} q^{-i} (\bar \xi_n - [i])) = p^r \sqcap^{r-1}_{i=0} [n-i], 
  \tag{2.46} 
\end{equation}
a $q$-analog of the familiar formula
\begin{equation}
  < \xi_n (\xi_n - 1) ... (\xi_n - r+1) > = p^r n ... (n - r+1). \tag{2.47} 
\end{equation}

Second, let
\begin{equation}
  \mu^\prime_r = < \xi^r_n >, \ \ \mu_r = < (\xi_n - < \xi_n > )^r >, \ \ \ r \in \Z_+, 
  \tag{2.48} 
\end{equation}
be the moments, around zero and $< \xi_n>$ respectively, considered as
functions of $p, n, r$.  Romanovsky [9] has proved that 
\begin{gather}
  \mu^\prime_{r+1} = (n p + p (1 - p) {d \over dp}) (\mu^\prime_r), 
  \tag{2.49} \\
  \mu_{r+1} = p(1 - p) (n r \mu_{r-1} + {d \mu_r \over dp}). \tag{2.50} 
\end{gather}
Formula (2.49) has the following $q$-analog:
\begin{equation}
  \mu^\prime_{r+1} = ([n] p + p (1-p) {d \over d_q p} ) (\mu^\prime_r). \tag{2.51} 
\end{equation}
Formula (2.50) has no clear $q$-analog, in part because the notion of 
{\it{higher \ central \ moments}} is not unique in $q$-probability.  Certainly, the classical 
definition
\begin{equation}
  \mu_r = E ((\xi - < \xi >)^r ) \tag{2.52} 
\end{equation}
is not useful, as the objects
\begin{equation}
  {d \over d_q x} ((\alpha + \beta x)^r) \tag{2.53} 
\end{equation}
lie outside compact formulae of $q$-analysis.  Two other possible definitions
are
\begin{equation}
  \mu_r = E(( \ \xi \ \lower1.25ex\hbox{$\buildrel-\over{{q}}$} < \xi > )^r), \tag{2.54} 
\end{equation}
where [7]
\begin{equation}
  (a \lower1.25ex\hbox{$\buildrel+\over{{q}}$} b)^n = \sum^n_{\kappa =0} 
  \bigg[{n \atop \kappa}\bigg] 
  a^\kappa b^{n-\kappa}, \ \ \ n \in \Z_+, \tag{2.55} 
\end{equation}
and
\begin{equation}
  \mu_r (s) = E ((\xi -^{\hskip-.076truein \cdot} \ {{q^s<}\xi >)^r}). \tag{2.56} 
\end{equation}
Using the definition (2.56) for the $q$-Bernoulli distribution (2.4), we find
\begin{gather}
  \mu_{r+1} (0; p) = p(1-p) (q [n] [r] \mu_{r-1} (1; pq) + {d \over d_q p} (\mu_r (1 ; p))), 
  \tag{2.57} \\
  \mu_{r+1} (-r; p) = p(1-p) (q^{-r} [n] [r] \mu_{r-1} (- r ; pq) + {d \over
  d_q p} (\mu_r 
  (-r; p ))). 
  \tag{2.58} 
\end{gather}
These formulae indicate that the central limit theorem for the $q$-Bernoulli
distribution  may not exist, at least in the classical sense.

\section{${\pmq}$-analogs of negative binomial distributions}

The negative binomial distribution, also called Pascal distribution, can be arrived at via 
many different routes.  Perhaps the simplest one is as the waiting time in a succession of 
Bernoulli trials until the appearance of $r^{th}$ non-zero for the first time.

Let's first consider the case $r=1$.  Let $W$ be the random variable, waiting
time until first non-zero.  By formula (2.29), 
\begin{equation}
  Pr (W = j) = Pr ({\bf{0}}^{j-1} q^0) = (1 -^{\hskip-.076truein \cdot} p)^{j-1}
  q^{j-1} p, \ \ \ j \in \bn. \tag{3.1} 
\end{equation}
We further $q$-modify this formula by setting
\begin{equation}
  Pr (\bar W = [j]) = (1 -^{\hskip-.076truein \cdot} p)^{j-1} q^{j-1} p
  , \ \ \ j \in \bn. \tag{3.2} 
\end{equation}
This is our $q$-analog of the geometric distribution.  Since
\begin{equation}
  \sum^\infty_{j=1} (1 -^{\hskip-.076truein \cdot} p)^{j-1} q^{j-1} p = 1 - Pr 
  ({\bf{0}}^\infty) = 1 - (1 -^{\hskip-.076truein \cdot} p)^\infty, 
  \tag{3.3} 
\end{equation}
we get a $q$-analog of the formula for the sum of a geometric progression:
\begin{equation}
  \sum^\infty_{s=0} (1 -^{\hskip-.076truein \cdot} 
  p)^s q^s = {1 - (1 -^{\hskip-.076truein \cdot} p)^\infty \over p} . 
  \tag{3.4} 
\end{equation}
\textbf{Remark 3.5.} Most of the formulae appearing in this paper
remain true when 
$q$ is considered as a formal variable, or as a complex one (with occasional restrictions of 
the type $|q|<1, \ |q|>1$, etc.)  It is only for the sake of probability interpretations that 
$q$ is considered to be a real number between 0 and 1.

Formula (3.4) has a finite counterpart.
\begin{equation}
  \sum^N_{s=0} (1 -^{\hskip-.076truein \cdot} x)^s q^s = {1 -
  (1 -^{\hskip-.076truein \cdot} x)^{N+1} \over x}. \tag{3.5} 
\end{equation}
This relation is easily checked by induction on $N$.

For general $r \in \bn$, the probability that the $r^{th}$ non-zero occurs at
exactly the 
$j^{th}$ trial, $j \geq r$, is, by formulae (2.29), (2.30), (2.33): 
\begin{gather}
  \sum_{|a(r-1)|=j-r} Pr ({\bf{0}}^{a_{0}} q^0 ... {\bf{0}}^{a(r-1)} q^{r-1} )=
  \sum (1 -^{\hskip-.076truein \cdot}  p)^{j-r} 
  p^r q^{r(j-r)} q^{- \sum^{r-1}_{0} sa(s)}  \notag\\
  \qquad = (1 -^{\hskip-.076truein \cdot} p)^{j-r} p^r q^{r(j-r)}
  q^{-(j-r)(r-1)} \bigg[{j-r+r-1 \atop  r-1} \bigg]  \notag\\
  \qquad = (1 -^{\hskip-.076truein \cdot} p)^{j-r} q^{j-r} p^r 
  \bigg[{j-1 \atop r-1}\bigg]. \tag{3.6} 
\end{gather}
Since the probability of having exactly $\ell$ non-zeroes during an infinite
number of  trials is, by formula (2.21), 
\begin{equation}
  {1 \over [\ell]!} \bigg({p \over 1-q}\bigg) ^\ell (1 -^{\hskip-.076truein \cdot} 
  p)^\infty, \tag{3.7} 
\end{equation}
the probability of having $< r$ non-zeroes is, therefore, 
\begin{equation}
  \bigg(\sum^{r-1}_{\ell=1} {1 \over [\ell]!} \bigg({p \over 1 -q} \bigg)^\ell \bigg)
  (1 -^{\hskip-.076truein \cdot} p)^\infty. \tag{3.8} 
\end{equation}
Thus,
\begin{equation}
  \sum^\infty_{s=0} (1 -^{\hskip-.076truein \cdot} p)^s q^s \bigg[{r-1+s \atop s}\bigg] 
  = p^{-r} \bigg(1 - (1 -^{\hskip-.076truein \cdot} p)
  ^\infty \sum^{r-1}_{\ell=0} {1 \over [\ell]!} \bigg( {p \over 1-q}
  \bigg)^\ell \bigg). 
  \tag{3.9} 
\end{equation}
This identity can be gotten directly from formula (3.4) by applying the
operator 
${\dps{\bigg({d \over d_q p}\bigg)^{r-1}}}$ to it.

Notice that formulae (3.3) and (3.9) show that our $q$-distributions do not 
sum up to 1 and therefore have to be re-scaled.  For example, formula (3.2) becomes
\begin{equation}
  Pr (\bar W_q = [j]) = {1 \over 1-(1-^{\hskip-.076truein \cdot} p)^\infty} 
  p (1 -^{\hskip-.076truein \cdot} p)^{j-1} q^{j-1}, \ \ j-1 \in \Z_+. 
  \tag{3.10} 
\end{equation}
\textbf{Remark 3.11.} Formula (3.6) could have been arrived at via one of
the  standard routes as the conditional probability of having $r-1$ non-zeroes
at the $1^{st}$ 
$j-1$ throws, with probability ${\dps{\bigg[{j-1 \atop r-1}\bigg] p^{r-1} 
(1 -^{\hskip-.076truein \cdot} p)^{j-r}}} $ by formula (2.4), followed by a non-zero appearing 
at the $j^{th}$ throw, with probability $q^{j-r}p$ by formula (2.26b).

We conclude this section by re-visiting the ``probl$\grave{\rm{e}}$me de 
parties'', one of the first 
problems in probability discussed and solved by Fermat and Pascal in their correspondence.  
In essence, we want to find the probability that $a$ non-zeroes appear before $b$ zeroes in 
the Bernoulli trial of $a+b-1$ throws.  This can happen in either one of the $b$ ways, when 
the $a^{th}$ non-zero appears at the $(a+\ell)^{th}$ trial, $0 \leq \ell \leq b - 1$.  By 
formula (3.6), the probability of this is $(j = a + \ell, \ r = a)$:
\begin{equation}
  (1 -^{\hskip-.076truein \cdot} \ p)^\ell q^\ell p^a \bigg[{a + \ell -1 \atop a - 1}\bigg], 
  \tag{3.12} 
\end{equation}
so that the total probability is
\begin{equation}
  P_1 = p^a \sum^{b-1}_{\ell =0} \bigg[{a+\ell-1 \atop a-1}\bigg] (1-^{\hskip-.076truein \cdot}
  p)^\ell q^\ell. \tag{3.13} 
\end{equation}
Similarly, the event that $b$ zeroes appear before $a$ non-zeroes can happen
in one of the 
$a$ ways, when the $b^{th}$ zero appears at the $(b+\ell)^{th}$ trial, $0 \leq \ell \leq a-1$.  
The probability of this event is, by formula (2.26a):
\begin{equation}
  \bigg[{b+\ell-1 \atop b-1}\bigg] p^\ell (1 -^{\hskip-.076truein \cdot} p)^{b-1} (1 - q^{b-1} 
  p) = \bigg[{b+ \ell-1 \atop b-1}\bigg] p^\ell (1 -^{\hskip-.076truein \cdot}
  p)^b . 
  \tag{3.14} 
\end{equation}
Thus, the total probability of $b$ zeroes appearing before $a$ non-zeroes is 
\begin{equation}
  P_2 = (1 -^{\hskip-.076truein \cdot} p)^b \sum^{a-1}_{\ell=0} \bigg[{b+\ell-1 \atop b-1} \bigg]
  p^\ell. \tag{3.15} 
\end{equation}
Since $P_1 + P_2 = 1$, we find
\begin{equation}
  p^a \sum^{b-1}_{\ell=1} \bigg[{a+\ell-1 \atop a-1} \bigg] (1 -^{\hskip-.076truein \cdot} p)^
  \ell q^\ell + (1 -^{\hskip-.076truein \cdot} p)^b \sum^{a-1}_{\ell=0}
  \bigg[{b+\ell-1 
  \atop b-1} \bigg] p^\ell =1, \ \ \forall a, b \in \bn, \tag{3.16} 
\end{equation}
an identity which is not immediately obvious.

The same probability $P_1$ can be calculated differently, as the outcome, 
out of $a + b-1$ trials, 
of $a+s$ non-zeroes, $0\leq s \leq b-1$.  Thus,
\begin{equation}
  P_1 = \sum^{b-1}_{s=0} \bigg[{a+b-1 \atop a+s} \bigg] p^{a+s}
(1 -^{\hskip-.076truein \cdot}
  p)^{b-1-s}, \tag{3.17} 
\end{equation}
and we arrive at another nonobvious (even for $q=1)$ identity
\begin{equation}
  p^a \sum^{b-1}_{\ell=0} \bigg[{a+\ell-1 \atop \ell } \bigg]
(1-^{\hskip-.076truein \cdot}
  p)^\ell q^\ell = \sum^{b-1}_{s=0} \bigg[{a+b-1 \atop a+s}
  \bigg] p^{a+s} (1 -^{\hskip-.076truein \cdot}
  p)^{b-1-s}. \tag{3.18} 
\end{equation}

\section{{$\pmq$}-Poisson distribution}

One of the shortest derivations of the Poisson distribution consists of considering, as Poisson 
originally did, the limit
\begin{equation}
  n \rightarrow \infty, \ \ \ pn \rightarrow \lambda \tag{4.1} 
\end{equation}
in the Bernoulli distribution:
\begin{equation}
  Pr (\xi = \kappa) = {n \choose \kappa} p^\kappa (1 - p)^{n-k} = {n ... (n - k+1) \over 
  \kappa!} {\lambda^\kappa \over n^\kappa} (1 - {\lambda \over n})^{n-\kappa}
  \rightarrow  {\lambda^\kappa \over \kappa!}  e^{-\lambda} . \tag{4.2} 
\end{equation}
The $q$-picture is more interesting. First, for $|q|<1$, the expression 
$Pr (\bar \xi_n = [\kappa]) $ (2.4) has the $n \rightarrow \infty$ - limit (2.21): 
\begin{equation}
  \bigg[{n \atop \kappa}\bigg] p^\kappa (1 -^{\hskip-.076truein \cdot} p)^{n-\kappa} 
  \longrightarrow 
  {1 \over [k]!} \bigg({p \over 1-q} \bigg)^\kappa (1 -^{\hskip-.076truein
  \cdot} p)^\infty. 
  \tag{4.3} 
\end{equation}
We can get a $q$-Poisson distribution from this by setting
\begin{equation}
  p = \lim_{n\to\infty} {\lambda \over [n] } = \lambda (1-q), \tag{4.4} 
\end{equation}
so that
\begin{equation}
  Pr(X = [\kappa]) = {\lambda^\kappa \over [\kappa]!} (1 -^{\hskip-.076truein \cdot}
  \lambda (1 - q))^\infty, \tag{4.5} 
\end{equation}
and therefore
\begin{equation}
  E_0 (\lambda)= E_0 (\lambda ; q) = \sum^\infty_{\kappa = 0} {\lambda^\kappa \over [\kappa]!} 
  = {1 \over (1 -^{\hskip-.076truein \cdot} \lambda (1-q))^\infty}, \tag{4.6} 
\end{equation}
the well-known formula.  Here
\begin{equation}
  E_\mu (\lambda) = \sum^\infty_{\kappa =0} {\lambda^\kappa \over [\kappa]!} q^{\mu{\kappa 
  \choose 2}} \tag{4.7} 
\end{equation}
is the $q$-family of exponentials:
\begin{equation}
  {d E_\mu (\lambda) \over d_q \lambda} = E_\mu (q^\mu \lambda), \ \ \ E_\mu (0) = 1. 
  \tag{4.8} 
\end{equation}
By Euler's formula (2.11), 
\begin{equation}
  (1 -^{\hskip-.076truein \cdot} p)^\infty = \sum^\infty_{j=0} \bigg(- {p \over 1-q}\bigg)^j 
  {1 \over [j]!} q^{{j \choose 2}} = E_1 \bigg(- {p \over 1-q} \bigg), 
  \tag{4.9} 
\end{equation}
so that
\begin{equation}
  E_0 (\lambda) E_1 (- \lambda) = 1; \tag{4.10} 
\end{equation}
the latter 2 formulae are of course classic.

Formula (4.10) can be generalized, as follows.  Consider the probability
generating function  for the $q$-Bernoulli distribution:
\begin{gather}
  F_{B;n} (z) = \sum^\infty_{k=0} z^\kappa Pr (\bar \xi_n = [\kappa]) = \sum^n_{\kappa =0} 
  \bigg[{n \atop \kappa}\bigg] z^\kappa p^\kappa (1 -^{\hskip-.076truein \cdot}
  p)^{n-\kappa}  \notag\\
  \qquad=\sum \bigg[{n \atop \kappa}\bigg] (zp)^\kappa (1 -^{\hskip-.076truein
  \cdot} p)^{n-\kappa} 
  \ \ [{\rm{by}} \ (2.8)] \ 
  = \sum \bigg[{n \atop \kappa} \bigg] (z p
  -^{\hskip-.076truein \cdot}  p)^\kappa  \notag \\
  \qquad = \sum \bigg[{n \atop \kappa} \bigg]
  p^\kappa (z -^{\hskip-.076truein \cdot}1)^
  \kappa.\tag{4.11} 
\end{gather}
As $n \rightarrow \infty$, this generating function becomes
\begin{equation}
  F_{B;\infty} = \sum^\infty_{\kappa = 0} {1 \over [\kappa]!} \bigg( {p \over 1 - q} \bigg)^
  \kappa (z -^{\hskip-.076truein \cdot} 1)^\kappa. \tag{4.12} 
\end{equation}
On the other hand, since
\begin{equation}
  (1 -^{\hskip-.076truein \cdot} x)^{\alpha + \beta} = (1 -^{\hskip-.076truein \cdot} x)
  ^\alpha (1 -^{\hskip-.076truein \cdot} q^\alpha x)^\beta, \tag{4.13a} 
\end{equation}
we have:
\begin{equation}
  (1 -^{\hskip-.076truein \cdot} x)^{- \beta} = {1 \over (1 -^{\hskip-.076truein \cdot} q^
  {- \beta} x)^\beta}, \tag{4.13b} 
\end{equation}
and therefore
\begin{gather*}
  (1 -^{\hskip-.076truein \cdot} p)^{n-\kappa} = (1  -^{\hskip-.076truein \cdot} p)^n 
  (1 -^{\hskip-.076truein \cdot} q^n p)^{-\kappa} = {(1 -^{\hskip-.076truein
  \cdot} p)^n \over  (1 -^{\hskip-.076truein \cdot} q^{n-\kappa} p)^\kappa} \\
  \qquad\lower.25ex\hbox{$\buildrel n \to\infty\over{{\rightarrow}}$} \ \ 
  (1 -^{\hskip-.076truein \cdot} p)^\infty = E_1 \bigg( - {p \over 1 - q} \bigg)
  = {1 \over  E_0 (\lambda)}, 
\end{gather*}
so that
\begin{equation}
  F_{B; \infty} = \lim_{n \to \infty} \sum z^\kappa \bigg[{n \atop \kappa}\bigg] p^\kappa
  (1 -^{\hskip-.076truein \cdot} p)^{n-\kappa} = \sum \bigg({zp \over
  1-q}\bigg)^\kappa  {1 \over [\kappa]!} (1 -^{\hskip-.076truein \cdot}
  p)^\infty = {E_0 (\lambda z) \over E_0  (\lambda)}. \tag{4.14} 
\end{equation}
Thus, 
\begin{equation}
  {E_0 (\lambda z) \over E_0 (z)} = \sum^\infty_{\kappa = 0} {\lambda^\kappa \over 
  [\kappa]!} (z -^{\hskip-.076truein \cdot} 1)^\kappa, \tag{4.15} 
\end{equation}
which can be equivalently rewritten as
\begin{equation}
  E_0 (b) E_1 (-a) = \sum^\infty_{\kappa=0} {(b -^{\hskip-.076truein \cdot} a)^\kappa 
  \over [\kappa]!}. \tag{4.16} 
\end{equation}

Taking the limit $n \rightarrow \infty$ in formulae (2.17), (2.20), we find:
\begin{gather}
  < X > = \lambda, \tag{4.17} \\
  Var (X) = \lambda (1 - (1 - q) \lambda). \tag{4.18} 
\end{gather}

Let us now consider the case when $|q| > 1$, so that we are walking outside the traditional 
probability theories.  Again,  set
\begin{equation}
  p = {\lambda \over [n]} . \tag{4.19} 
\end{equation}
Then, as $n \rightarrow \infty$, 
\begin{equation}
  \bigg[{n \atop \kappa}\bigg] p^\kappa = {[n] ... [n-\kappa+1] \over [\kappa]!} {\lambda^
  \kappa \over [n]^\kappa} \rightarrow {\lambda^\kappa \over [\kappa]!}
  q^{-{\kappa \choose 2}}, 
  \tag{4.20} 
\end{equation}
since, as $n \rightarrow \infty,$
\begin{equation}
  {[n - \ell] \over [n]} = {q^{n-\ell} -1 \over q^n -1} \rightarrow q^{- \ell}. \tag{4.21} 
\end{equation}
Next, by formula (4.13b), 
\begin{equation}
  (1 -^{\hskip-.076truein \cdot} p)^{n-\kappa} = (1 -^{\hskip-.076truein \cdot} {\lambda \over 
  [n]})^{n-\kappa} = {(1 -^{\hskip-.076truein \cdot} {\lambda \over [n]})^n
  \over  (1  -^{\hskip-.076truein \cdot} q^{n-\kappa} {\lambda \over
  [n]})^\kappa} \rightarrow 
  \lim_{n\to\infty} (1 -^{\hskip-.076truein \cdot} {\lambda \over [n]})^n (1
  -^{\hskip-.076truein \cdot}
  \lambda)^{- \kappa}. \tag{4.22} 
\end{equation}
Thus,
\begin{equation}
  Pr (X = [\kappa]) = {q^{-{\kappa \choose 2}} \over [\kappa]!} \lambda^\kappa (1 
  -^{\hskip-.076truein \cdot} \lambda)^{-\kappa} \lim_{n \rightarrow \infty} (1
  -^{\hskip-.076truein \cdot}
  {\lambda \over [n]})^n, \ \ |q|>1. \tag{4.23} 
\end{equation}
In particular,
\begin{equation}
  \lim_{n \to\infty} (1 -^{\hskip-.076truein \cdot} {\lambda \over [n]})^n = 
  \bigg(\sum^\infty_{\kappa = 0} {q^{-{\kappa \choose 2}} \over [\kappa]!}
  \lambda^\kappa  (1 -^{\hskip-.076truein \cdot} \lambda)^{-\kappa} \bigg)^{-1}
  , \ \ \ |q| > 1; \tag{4.24} 
\end{equation}
this formula is not a $q$-analog of anything classical.

\section{${\pmq}$-hypergeometric distribution}

Imagine that we have an urn consisting of two types of balls: $m$ marked `1' and $u$ marked 
`0'.  We pick out at random one ball, record its value and {\it{leave \ it \ outside \ the}} \
 $urn$; 
then proceed again, for a total of $n$ draws.  Had we returned each picked ball back into the 
urn, we would have the Bernoulli trials; since we don't return the balls, we get something 
different, called the hypergeometric distribution:  the probability of ending up with 
$\kappa$ `1' balls out of $n$ draws is
\begin{equation}
  Pr (\xi_n = \kappa) = \bigg({m \atop \kappa}\bigg) \bigg({u \atop n-\kappa}\bigg)\bigg/
  \bigg({N \atop n}\bigg), \ \ \ N = m+u, \tag{5.1} 
\end{equation}
see [4, 5].

As a $q$-analogue of this distribution we set
\begin{equation}
  Pr (\bar \xi_n = [\kappa]) = \bigg[{m \atop \kappa}\bigg] \bigg[{u \atop n-\kappa}\bigg] 
  q^{(m - \kappa) (n - \kappa)} \bigg/ \bigg[{N \atop n } \bigg]. \tag{5.2} 
\end{equation}
To justify this definition we have to verify that 
\begin{equation}
  \sum_\kappa \bigg[{m \atop \kappa}\bigg] \bigg[{u \atop n-\kappa}\bigg] q^{(m-\kappa)(n-
  \kappa)} = \bigg[{m+u \atop n }\bigg]. \tag{5.3} 
\end{equation}
This identity results by picking the $x^n$-coefficient in formula (4.13a):
\begin{equation}
  (1 -^{\hskip-.076truein \cdot} x)^m (1 -^{\hskip-.076truein \cdot} q^m x)^u = 
  (1 -^{\hskip-.076truein \cdot} x)^{m+u}, \tag{5.4} 
\end{equation}
and using the Euler formula (2.11):
\begin{gather}
  \sum_\kappa \bigg[{m \atop \kappa} \bigg] (- x)^\kappa q^{{\kappa \choose 2}} \sum_\ell 
  \bigg[{u \atop \ell}\bigg] (- x)^\ell q^{m \ell} q^{{\ell \choose 2}}
  \notag\\
  \qquad = \sum_n  (-x)^n q^{{n \choose 2}} \sum_{k + \ell =n} \bigg[{m \atop
  \kappa}\bigg] \bigg[{u \atop 
  \ell}\bigg] q^{{\kappa \choose 2} + {\ell \choose 2} - {\kappa + \ell \choose
  2}} q^{m\ell} \notag\\
  \qquad= \sum_n (-x)^n q^{{n \choose 2}} \sum_\kappa \bigg[{m \atop
  \kappa}\bigg]
  \bigg[  {u \atop n-\kappa}\bigg] q^{(m-\kappa)(n-\kappa)} = \sum_n \bigg[{m +
  u \atop n}\bigg]  (-x)^n q^{{n \choose 2}}, \tag{5.5} 
\end{gather}
where we used the obvious relation
\begin{equation}
  \bigg({\kappa + \ell \atop 2}\bigg) = \bigg({\kappa \atop 2}\bigg) + \bigg({\ell \atop 2} 
  \bigg) + \kappa \ell. \tag{5.6} 
\end{equation}

Similar to the Bernoulli case, we can treat the $q$-hypergeometric
distribution (5.2) as a  macroscopic object and inquire about its microscopic
representation.  The latter can be  guessed from the relations
\begin{gather}
  P (\bar \xi_n = [n]) = \bigg[{m \atop n}\bigg] \bigg/ \bigg[{N \atop n}\bigg] = {[m] 
  \over [n]} {[m -1] \over [N-1]} ... {[m-n+1] \over [N-n+1]} , \tag{5.7a} \\
  Pr (\bar \xi_n = 0) = \bigg[{u \atop n}\bigg] q^{mn}/\bigg[{N \atop n}\bigg]
  = 
  {[u] \over [N]} q^m {[u-1] \over [N-1]} q^m ... {[u-n+1] \over [N-n+1]} q^m,
  \tag{5.7b} 
\end{gather}
which suggest that in the representation
\begin{equation}
  \bar \xi_n = \bar \zeta_1 + ... + \bar \zeta_n \tag{5.8} 
\end{equation}
of $n$ successive draws, we should set
\begin{gather}
  Pr ({\bf{0}}^{a(0)} q^0... {\bf{0}}^{a(k-1)} q^{k-1} {\bf{0}}^{a(\kappa)})=q ^{\sum_{0}^{\kappa}
  (m-s)a(s)} \bigg[{m \atop \kappa}\bigg] \bigg[{u \atop n-\kappa}\bigg]
  \bigg/ 
  \bigg[{N \atop n}\bigg] \bigg[{n \atop \kappa}\bigg] \tag{5.9a} \\
  \qquad= q^{\sum^\kappa_{0}(m-s)a(s)} \sqcap^{\kappa-1}_{i=0} [m-i]
  \sqcap^{n-\kappa-1}_{j=0}  [u-j] / \sqcap^{n-1}_{\gamma =0} [N-\gamma]. 
  \tag{5.9b}  \\
  Pr({\bf{0}}^a) = q^{ma} \bigg[{u \atop a}\bigg] \bigg/ \bigg[{N \atop
  a}\bigg]. \tag{5.9c} 
\end{gather}
(In the $\bar \zeta$-language, we have
\begin{equation}
  Pr (\bar \zeta \not=0) = {[m] \over [N]} , \ \ \ Pr (\bar \zeta = 0) = q^m {[u] \over [N] }
  \tag{5.10} 
\end{equation}
at each ball pick-out when the number of marked `nonzero' balls in the urn is
$m$, and  the number of those marked `zero' is $u = N-m$.) 

To prove that microscopic formula (5.9a) implies the macroscopic formula
(5.2), we need  to verify that 
\begin{equation}
  \sum_{|a(\kappa)|=n-\kappa} q^{\sum_{0}^{\kappa} (m-s)a(s)} = q^{(m-\kappa)(n-\kappa)}
  \bigg[{n \atop \kappa}\bigg], \tag{5.11} 
\end{equation}
and this equality follows from the already proven formula (2.35):
\begin{align*}
  \sum_{|a(\kappa)|=n-\kappa} q^{\sum^{\kappa}_{0} (m-s)a(s)} &= \sum q^{m\sum
  a(s)} 
  q^{-\sum sa(s)} = q^{m(n-\kappa)} q^{- (n-\kappa)\kappa} 
  \bigg[{n - \kappa + \kappa \atop \kappa}\bigg] \\
  &= q^{(m-\kappa)(n-\kappa)} \bigg[{n \atop \kappa}\bigg].
\end{align*}

In the classical case $q =1$, we can rewrite formula (5.1) as 
\begin{equation}
  Pr (\xi_n = \kappa) = \bigg({n \atop \kappa}\bigg) \sqcap^{k-1}_{s=0} {m-s \over N-s} 
  \sqcap^{n-\kappa+1}_{\ell=0} {u - \ell \over N - \kappa -\ell}, \tag{5.12} 
\end{equation}
so that in the limit
\begin{equation}
  N \rightarrow \infty, \ \ {m \over N} \rightarrow p, \ \ {u \over N} \rightarrow 1-p, 
  \tag{5.13} 
\end{equation}
formula (5.12) becomes the Bernoulli one; the classical explanation is that
when $m$ and 
$u$ are  large, it makes little difference whether the picked-out balls are returned back to the 
urn or not.

For $q \not= 1$, the situation is more interesting.  Certainly formulae (5.13)
are not the  correct ones.  We proceed as follows.

Let $|q|>1$.  Set
\begin{equation}
  N \rightarrow \infty, \ \ {[m] \over [N]} \rightarrow p, \tag{5.14} 
\end{equation}
and re-write formula (5.2) in the form
\begin{equation}
  Pr (\bar \xi_n = [\kappa])= \bigg[{n \atop \kappa}\bigg] \sqcap^{\kappa-1}_{s=0} 
  \bigg({[m-s] q^{\kappa -n} \over [N - (n-\kappa)-s]} \bigg)
  \sqcap^{n-\kappa-1}_{\ell=0} 
  \bigg({[u - \ell]q^m \over [N-\ell]}\bigg). \tag{5.15} 
\end{equation}
Now,
\begin{gather}
  {[m-s] q^{k-n} \over [N - (n-\kappa)-s]} = q^{\kappa -n} {[-s] + q^{-s} [m] \over 
  [\kappa -s-n] + q^{\kappa -s-n} [N]} \ \ {\rm{[by \ (5.14)]}} \ 
  \notag \\
  \qquad \rightarrow 
  q^{\kappa -n} {q^{-s} \over q^{\kappa -s-n}} p = p, 
  \tag{5.16a} \\
  {[u-\ell[ q^m \over [N - \ell]} = {[N-m-\ell]q^m \over [N-\ell]} =
  {[N-\ell]- [m] \over 
  [N-\ell]} = 1  - {[m] \over [-\ell]+q^{-\ell}[N]} \rightarrow 1-q^\ell p, 
  \tag{5.16b} 
\end{gather}
so that
\begin{equation}
  Pr (\bar \xi_n = [\kappa]) \rightarrow \bigg[{n \atop \kappa}\bigg] p^\kappa (1 
  -^{\hskip-.076truein \cdot} p)^{n-\kappa}, \tag{5.17} 
\end{equation}
as desired.

Suppose now that $|q|<1$.  Denote that $q$ by $Q$.  Set 
\begin{equation}
  q=Q^{-1}, \ \ \ |q|>1. \tag{5.18} 
\end{equation}
Since
\begin{gather}
  [x]_Q = q^{1-x} [x]_q, \tag{5.19a} \\
  [k]!_Q = q^{-{\kappa \choose 2}} [\kappa]!_q, \tag{5.19b} \\
  \bigg[{a \atop b}\bigg]_Q = q^{b(b-a) } \bigg[{a \atop b}\bigg]_q, 
  \tag{5.19c} 
\end{gather}
we have
\begin{align}
  Pr (\bar \xi_n &= [\kappa]_q) = Q^{(m-\kappa)(n-\kappa)} \bigg[{m \atop
  \kappa}\bigg]_Q 
  \bigg[{u \atop n-\kappa}\bigg]_Q\bigg/\bigg[{N \atop n}\bigg]_Q \notag \\
  &= q^{(\kappa-\kappa^\prime)(n-\kappa^\prime)} \bigg[{u \atop
  \kappa^\prime}\bigg] _q
  \bigg[{m \atop n-\kappa^\prime}\bigg]_q\bigg/ \bigg[{N \atop n }\bigg]_q,
  \tag{5.20} 
\end{align}
where
\begin{equation}
  \kappa^\prime = n-\kappa. \tag{5.21} 
\end{equation}
We see that our formula (5.2) in the form 
\begin{equation}
  Pr (\tilde \xi_n = \kappa) = \bigg[{m \atop \kappa}\bigg] \bigg[{u \atop n - \kappa}\bigg] 
  q^{(m-\kappa)(n-\kappa)} \bigg/ \bigg[{N \atop n}\bigg] \tag{5.22} 
\end{equation}
allows the symmetry
\begin{equation}
  m \rightarrow u, \ u \rightarrow m, \ \kappa \rightarrow n-\kappa, \ q \rightarrow q^{-1}. 
  \tag{5.23} 
\end{equation}
Setting
\begin{equation}
  N \rightarrow \infty, \ {[u] \over [N]} \rightarrow p^\prime \ (=1 -p), \tag{5.24} 
\end{equation}
we find, in the same way as formula (5.17) was gotten, that
\begin{equation}
  Pr (\tilde \xi_n = [\kappa]) = \bigg[{n \atop \kappa}\bigg] (1 -^{\hskip-.076truein \cdot}
  p^\prime)^\kappa p^{\prime n-\kappa}, \tag{5.25} 
\end{equation}
a different $q$-version of the classical Bernoulli distribution.

\section{$\pmq$-contagious distribution}

Suppose we again, like in the preceding section, have an urn with $m$ marked
and $u$ unmarked  balls.  We pick out one ball at random, record its value,
and then return to the urn $s+1$  balls identical to the one we just picked
out.  If $s=0$, we return the ball itself, and this is  the Bernoulli scheme;
if $s=-1$, we return nothing, and this is the hypergeometric scheme.   For
general $s$, the probability to pick out $k$ marked (by `1') balls out of $n$
draws is

\begin{equation}
  Pr (\xi_n = \kappa) = \bigg({n \atop \kappa}\bigg) \sqcap^{\kappa-1}_{\alpha =0} 
  (m + as) \sqcap^{n-\kappa-1}_{\beta =0} (u+\beta s) \bigg/
  \sqcap^{n-1}_{\gamma=0}  (m+u+\gamma s). \tag{6.1} 
\end{equation}
This particular member of the family of the so-called ``contagious
distributions'' was  discovered by Eggenberger and P${\rm{\acute{o}}}$lya in
1923 [2, 8].

As a $q$-analog of this distribution, we set
\begin{gather}
  Pr (\bar \xi_n = [\kappa]) \notag\\
  \qquad= \bigg[{n \atop \kappa}\bigg]_{q^{-s}} {q^{(m+s
  \kappa) (n -\kappa)}}
  \sqcap^{\kappa -1}_{\alpha=0} [m + \alpha s] \sqcap^{n-\kappa-1}_{\beta =0}
  [u+\beta s] 
  \bigg/ \sqcap^{n-1}_{\gamma=0} [m+u+\gamma s]. \tag{6.2} 
\end{gather}
Obviously, for the case $s=-1$ we recover the hypergeometric formula (5.2). 
For the case $s=0$,  we recover the {\it{classical }} Bernoulli formula (2.3)
with $p = [m]/[m+u])$,  {\it{not}} the 
$q$-Bernoulli formula (2.4).

To justify formula (6.2), we need to check that
\begin{equation}
  \sum^n_{\kappa=0} \bigg[{n \atop \kappa}\bigg]_{q-s} q^{(m+s\kappa)(n-\kappa)} 
  \sqcap^{\kappa -1}_{\alpha =0} [m + \alpha s] \sqcap^{n-\kappa-1}_{\beta =0}
  [u+\beta s] = 
  \sqcap^{n-1}_{\gamma=0} [m+u +\gamma s]. \tag{6.3} 
\end{equation}
To do that, we assume that $s\not=0$ and start with the formula (5.3) in the base $Q=q^{-s}$:
\begin{equation}
  \sum_{\kappa=0}^n \bigg[{M \atop \kappa}\bigg]_Q \bigg[{U \atop n-\kappa}\bigg]_Q 
  Q^{(M-\kappa)(n-\kappa) } = \bigg[{M + U \atop n}\bigg]_Q, \ \ \ Q = q^{-s},
  \tag{6.4} 
\end{equation}
which we rewrite as
\begin{equation}
  \sum^n_{\kappa=0} \bigg[{n \atop \kappa}\bigg]_Q Q^{(M-\kappa)(n-\kappa)} 
  \sqcap^{\kappa-1}_{\alpha=0} [M - \alpha]_Q \sqcap^{n-\kappa-1}_{\beta =0} 
  [U- \beta] _Q = \sqcap^{n-1}_{\gamma=0} [M + U - \gamma]_Q. \tag{6.5} 
\end{equation}
Multiplying both sides by $([ - r]_q)^n$, setting
\begin{equation}
  M = -  m/s, \ \ U = - u/s, \tag{6.6} 
\end{equation}
and noticing that 
\begin{equation}
  Q^{(M-\kappa)(n-\kappa)} = q^{-s(-\kappa - m/s)(n - \kappa)} = q^{(m+\kappa s) (n-\kappa)}, 
  \tag{6.7} 
\end{equation}
we arrive at formula (6.3).  The latter formula can be considered as a new
$q$-analog of  Newton's binomial.  

Similar to the hypergeometric case, we can arrive at the $q$-contagious
distribution (6.2)  as a macroscopic object, 
\begin{equation}
  \bar\xi_n = \bar\zeta_1 + ... + \bar \zeta_n \tag{6.8} 
\end{equation}
from the microscopic formulae
\begin{gather}
  Pr ({\bf{0}}^{a(0)} q^0 ... {\bf{0}}^{a(k-1)} q^{k-1} {\bf{0}}^{a(k)} ) 
  \notag \\
  \qquad = q^{\sum^{\kappa}_{i=0} (m+si)a(i)} \sqcap^{k-1}_{\alpha=0} [m+\alpha
  s]\sqcap^{n-\kappa-1}_{
  \beta=0} [u+\beta s] \bigg/ \sqcap^{n-1}_{\gamma=0} [m+u+\gamma s] ; 
  \tag{6.9} 
\end{gather}
the latter formulae are suggested by the extreme cases $k = n$ and $k=0$ of
formula (6.2):
\begin{gather}
  Pr (\bar \xi_n = [n])=\sqcap^{n-1}_{\alpha=0} 
  {[m+\alpha s] \over [N + \alpha s]}, 
  \ \ \ N = u + m, \tag{6.10a} \\
  Pr (\bar \xi_n = 0) = \sqcap^{n-1}_{\beta=0} 
  \bigg({[u + \beta s] \over [N + \beta s] } q^m \bigg). \tag{6.10b}
\end{gather}
To verify that microscopic formulae (6.9) imply the macroscopic formula (6.2),
we need to  check that 
\begin{equation}
  \sum_{|a(\kappa)|=n-\kappa} q^{\sum^{\kappa}_{i=0} (m+si)a(i)} = q^{(m + s \kappa)(n-\kappa}) 
  \bigg[{n \atop \kappa}\bigg]_{q^{-s}}. \tag{6.11} 
\end{equation}
Now, for the LHS of formula (6.11) we get:
\begin{gather}
  \sum_{|a(\kappa)|=n-\kappa} q^{\sum^{\kappa}_{i=0} (m+si)a(i)} 
  = q^{m (n-\kappa)} 
  \sum (q^{-s})^{-\sum i a(i)} \ \ {\rm{[by \ (2.35)}}] \notag\\
  \qquad = q^{m (n-\kappa)} (q^{-s})^{-(n-\kappa)\kappa} 
  \bigg[{n \atop \kappa}\bigg]_{q^{-s}} = q^{(m+s \kappa) (n-\kappa) } 
  \bigg[{n \atop \kappa}\bigg]_{q^{-s}}, \tag{6.12} 
\end{gather}
and this is exactly the RHS of formula (6.11).

\section{${\pmq}$-uniform distribution}

A classical random variable $X$ taking $M+1$ discreet values $v_0 < ... < v_M$, 
 each with equal probability 
$1/(M+1)$, represents a discrete uniform distribution.  The values $v_i'$s are
immaterial  and can be taken as $v_i = i, $ or $v_i = a+ hi$, or $v_i = [i]$,
\ ...

As a $q$-analog of this distribution, we set
\begin{equation}
  Pr (\tilde X = i) = q^i/[M+1], \ \ 0 \leq i \leq M, \tag{7.1} 
\end{equation}
or
\begin{equation}
  Pr (\bar X = [i]) = q^i/[M+1], \ \ \ 0 \leq i \leq M, \tag{7.2} 
\end{equation}
so that
\begin{gather}
  <\bar X> = {q [M] \over [2]}, \tag{7.3} \\
  <\bar X^2> = {q[M] [M+1] (q[2][M]+1) \over [2] [3]}, \tag{7.4} \\
  Var (\bar X) = {q[M] (q^2 [M] + [2]) \over [2]^2 [3]}. \tag{7.5} 
\end{gather}

Consider $n$ independent identically distributed, via the discrete uniform
distribution, random  variables $X_1, ..., X_n$.  The {\it{range}} of these
variables is the quantity
\begin{equation}
  r_n = \lower1.75ex\hbox{$\buildrel max\over{{i}}$}   
  (X_i) -  \lower1.75ex\hbox{$\buildrel min\over{{i}}$}  
  (X_i), \ \ \ 0 \leq r \leq M. \tag{7.6} 
\end{equation}
The random variable $r_n$ has the following distribution ([4]) p. 240):
\begin{gather}
  Pr (r_n=0) = 1/(M+1)^{n-1}, \tag{7.7a} \\
  Pr (r_n=\ell) = ((\ell+1)^n - 2 \ell^n + (\ell-1)^n) (M+1-\ell)/(M+1)^n, \ \
  \ 1 \leq \ell 
  \leq M. \tag{7.7b} 
\end{gather}
This distribution is certainly different from those appearing in the preceding
sections.

Let us calculate the $q$-analog of the distribution (7.7).  Taking as our
basic definition  formula (7.1), we have:
\begin{equation}
  Pr (r_n=0) = \sum^M_{i =0} (Pr (\tilde X = i))^n = 
  [M+1]_{{\dps{q^{n}}}}\bigg/ [M+1]^n ; \tag{7.8} 
\end{equation}
\begin{gather}
  Pr (r_n=1) = \sum^{M-1}_{i=0} \sum_{k \not= 0,n} \bigg({n \atop \kappa }\bigg) \bigg({
  q^i \over [M+1]}\bigg)^\kappa
  \bigg({q^{i+1} \over [M+1] }\bigg)^{n-\kappa} \notag \\
  \qquad = {1 \over [M+1]^n} \sum^{M-1}_{i=0} \ \ (\sum^n_{k=0} 
  \bigg({n \atop \kappa}\bigg) (q^i)^\kappa (q^{i+1})^{n-\kappa} - 
  (q^{i+1})^n - (q^i)^n) \notag \\
  \qquad= {1 \over [M+1]^n} \sum^{M-1}_{i = 0} ((q^i + q^{i+1} 
  )^n - q^{in} q^n - q^i) \notag \\
  \qquad= {1 \over [M+1]^n} \sum^{M-1}_{i = 0} (q^i)^n ([2]^n -
  [2]_{{\dps{q^{n}}}}) = {[M]_{\dps{q^{n}}} ([2]^n - [2]_{\dps{q^{n}}}) 
  \over [M+1]^n }; \tag{7.9} 
\end{gather}
finally, for $\ell \geq 2$, 
\begin{gather}
  Pr (r_n = \ell) = \sum^{M-\ell}_{i=0} \sum_{\kappa(0), \kappa(\ell)\not=0}
  {n! \over \kappa (0)!... \kappa 
  (\ell)!} \bigg({q^i \over [M+1]}\bigg)^{k(0)} ... \bigg({q^{i+\ell} \over
  [M+1]} \bigg)^ {\kappa (\ell)} \notag \\
  \qquad = {1 \over [M+1]^n} \sum^{M-\ell}_{i = 0} \bigg(\sum_{{\rm{all}} \
  \kappa' s}
  - \sum_{\kappa(0)=0} - \sum_{
  \kappa(\ell)=0} + \sum_{\kappa(0)=\kappa(\ell)=0} \bigg) \notag \\
  \qquad = {1 \over [M+1]^n}
  \sum^{M-\ell}_{i=0} \begin{array}[t]{l}
    \big((q^i +...+ q^{i+\ell})^n - (q^{i+1} + ... + q^{i+\ell})^n \\
    \ {}- (q^i + ... + q^{i+\ell-1})^n + (q^{i+1} + ... +
    q^{i+\ell-1})^n\big) 
  \end{array}\notag \\
  \qquad = {1 \over [M+1]^n} \sum^{M-\ell}_{i=0} q^{in} \bigg([\ell+1]^n - q^n
  [\ell]^n -  [\ell]^n + q^n [\ell-1]^n\bigg) \notag \\
  \qquad = {[M+1-\ell] _{\dps{q^{n}}} \over [M+1]^n} \bigg([\ell + 1]^n -
  [2]_{\dps{q^{n}}}
  [\ell]^n + q^n [\ell -1]^n 
  \bigg). \tag{7.10} 
\end{gather}
For $\ell=1$, formula (7.10) reproduces formula (7.8).  Thus, formulae
\begin{gather}
  Pr (r_n=0) = {[M+1]_{\dps{q^{n}}} \over [M+1]^n}, \tag{7.11a} \\
  Pr (r_n =\ell) = {[M+1-\ell]_{\dps{q^{n}}} \over [M+1]^n} 
  \bigg([\ell +1]^n - [2]_{\dps{q^{n}}}
  [\ell]^n + q^n [\ell-1]^n \bigg), \ \ 1 \leq \ell \leq M, \tag{7.11b} 
\end{gather}
are $q$-analogs of formulae (7.7).  Notice, that
\begin{equation}
  Pr (r_1 > 0) = 0, \ \ \ \forall q. \tag{7.12} 
\end{equation}

Although not immediately apparent, formulae (7.11) are not the {\it{only}}
$q$-analogs of  the classical formulae (7.7).  For example, for $n=2$, formulae
\begin{gather}
  Pr (r_2 = 0) = {1 \over [M+1]}, \tag{7.13a} \\
  Pr (r_2 = \ell) = {[2] \over [M+1]} (1 - {[\ell] \over [M+1]} )
  q^{M+1-2\ell} , 
  \ \ 1 \leq \ell \leq M, 
  \tag{7.13b} 
\end{gather}
are different from formulae $(7.11)\bigg|_{n=2}$.

\section{Concluding remarks}

The approach to $q$-Probability taken in this paper leaves the {\it{rules}} of
classical probability  intact and only $q$-deforms some basic probability
distributions.  It is quite likely that one  can develop some new/bizarre
rules of $q$-probability which contradict the comfortably  familiar intuition,
similar to what Quantum-mechanical interpretations  appear to a
Classical-mechanical disciple.  The most direct route to such new rules
probably goes through basic  {\it{continuous}} probability distributions when
integral $\int(\cdot)dx$ is replaced by the 
$q$-integral $\int (\cdot) d_qx$.

\def\uline{\textbf}

\label{kuper_lp}

\end{document}